\documentclass[11pt,a4paper]{amsart}
\usepackage[marginparwidth=0pt,margin=20truemm]{geometry} % margin

\usepackage{mypackage} % My packages.
\usepackage{mycommand} % My commands.

\usepackage{amscd}
\usepackage{url}
\usepackage{physics} % \dv{z} でd/dzを出すなどができるようになる

\DeclareMathOperator{\Rec}{Rec}
\DeclareMathOperator{\LHQX}{LHQX}
\DeclareMathOperator{\LHQW}{LHQW}
\DeclareMathOperator{\LHQWadm}{LHQW^*}
\DeclareMathOperator{\LHQE}{LHQE}
\DeclareMathOperator{\seg}{seg}
\DeclareMathOperator{\proc}{proc}
\DeclareMathOperator{\pav}{pav}

\usepackage{ctable} 
  {\end{enumerate}}

\numberwithin{equation}{section} % Setting of equation numbers 
\usepackage{mytheoremeng} % My theorem environments (English) with cleveref package

\usepackage{tikz}
\usetikzlibrary{arrows.meta, decorations, decorations.markings}
\begin{document}
% --------------------------------------------------------------------------

\title[Rectangulations avoiding a pattern]{Rectangulations avoiding a pattern}

\author[K. Sano]{Kaoru Sano}
\address{NTT Institute for Fundamental Mathematics, NTT Communication Science Laboratories, NTT, Inc., 2-4, Hikaridai, Seika-cho, Soraku-gun, Kyoto 619-0237, Japan}
\email{kaoru.sano@ntt.com}

\date{\today}

% --------------------------------------------------------------------------

\begin{abstract}
    Fix a strong rectangulation pattern $P$ of size $L$.
    We show that the growth constant of the class of strong rectangulations avoiding $P$ is strictly smaller than $\Lambda =27/2$, the growth constant for all strong rectangulations. More precisely, forbidding any such $P$ yields a pattern-uniform exponential drop of at least $\Lambda - 1/\Lambda^{3L-1}$.
    Consequently, the proportion of $P$-avoiding rectangulations among all rectangulations tends to zero as $n\to \infty$.
    This is the first result on the uniform drop of exponential growth for pattern-avoiding rectangulations.
    The proof utilizes the standard correspondence with leftmost history quadrant walks, along with a pattern-insertion scheme that controls the radius of convergence of the associated generating functions, thereby establishing the first uniform exponential upper bound for rectangulation classes defined by geometric avoidance.
\end{abstract}

\maketitle

% --------------------------------------------------------------------------

\tableofcontents

% --------------------------------------------------------------------------

% --------------------------------------------------------------------------

\section{Introduction} \label{sec: intro}

% --------------------------------------------------------------------------
A rectangulation is a partition of a rectangle into a finite number of rectangular blocks.
An $n$-rectangulation is a rectangulation consists of $n$ rectangles,
Two widely used equivalence relations lead to two counting problems:

In the weak sense, one identifies rectangulations that can be obtained from each other by sliding interior segment endpoints while preserving segment-rectangle incidence.
Let $\Rec^{\rm w}(n)$ be the set of equivalence classes of $n$-rectangulations in the weak sense.
A classical bijection due to \cite{ABP06} identifies weak rectangulations with Baxter permutations, and the classical study and recurrences for Baxter numbers by \cite{CGHK78} imply that $\# \Rec^{\rm w}(n)$ asymptotically behaves as
\[
    \# \Rec^{\rm w}(n) \sim \frac{32}{\pi \sqrt{3}} \frac{8^n}{n^4}\quad \text{as } n\to \infty.
\]
See \cite{OEIS_A001181} for the exact numbers of Baxter permutations.
Moreover, a bijection between Baxter permutations and bipolar orientations is given in \cite{BBMF08}.

In the strong sense, one identifies rectangulations up to homeomorphisms preserving the adjacencies of rectangles, disallowing cross junctions.
Let $\Rec^{\rm s}(n)$, or simply $\Rec(n)$, be the set of equivalence classes of $n$-rectangulations in the strong sense.
\cite{Rea12} exhibits a bijection between strong rectangulations and a pattern-defined permutation class, the $2$-clumped permutations, providing a clean combinatorial encoding.
This bijection is based on the study of \cite{BBMF08}.
The first systematic enumeration of generic rectangulations was performed by \cite{CM14}, which derived recursive formulas and computed the numbers of strong $n$-rectangulations for $n\le 28$.
Their numerical data suggested an exponential growth close to $(27/2)^n$.
\cite{Fus09} defined a transversal structure as a pair of non-intersecting bipolar orientations of a planar triangulation whose edges are colored and oriented according to two transversal directions (horizontal and vertical), and \cite{FNS24} calculated the asymptotics of the number of transversal structures.
Strong rectangulations correspond to transversal structures, which admit a powerful walk encoding, the leftmost history quadrant walk.
We recall the definition of leftmost history quadrant walks in \cref{subsec: notation}.
These bijections are now part of the standard dictionary.
A recent survey \cite{ACFF24} unifies these mappings and fixes notation.
By this correspondence and the result of \cite{FNS24}, the following sharp asymptotic is obtained.
\begin{thm}[{\cite[Theorem 4.2]{FNS24}}]\label{thm: asymptotics of rec n}
    There is a positive constant $c$ such that
    \begin{equation}\label{eq: asymptotics of n-rectangulations}
        \#\Rec(n) \sim c\Lambda^n n^{-\alpha} \quad\text{as }n\to \infty,
    \end{equation}
    where $\Lambda = 27/2$ and $\alpha = 1+\pi/\arccos(7/8)$.
\end{thm}

A {\it configuration of segments} is a finite set of horizontal and vertical line segments in the plane such that distinct segments neither cross nor share endpoints.
Two configurations $C$ and $C'$ are {\it weakly equivalent} if there exists a bijection $\gamma\colon C \longrightarrow C'$ satisfying the following conditions.
\begin{enumerate}
	\item $\gamma$ preserves orientation (horizontal/vertical).
	\item For any pair of segments, incidence is preserved, and the type of their contact is the same.
	\item Along each fixed segment, the order in which its neighbors appear is preserved on each side separately (left vs. right for vertical segments, and analogously for horizontal ones).
\end{enumerate}
They are {\it strongly equivalent} if they are weakly equivalent and satisfy the following additional condition.
\begin{enumerate}
    \item[(iv)]\label{item: pattern 4} for every segment, the single linear order of all its neighbors along the segment is preserved (i.e., the neighbors appear in the same order when traversing the segment from one end to the other, without distinguishing sides).
\end{enumerate}

A {\it weak} (resp. {\it strong}) {\it rectangulation pattern} is an equivalence class of configurations under weak (resp. strong) equivalence.
The {\it size} of a rectangulation pattern $P$ is the number of segments in any representative of $P$, and is denoted by $|P|$.
Given a rectangulation $R$, we represent it by the configuration of its interior maximal horizontal and vertical segments.
For a strong rectangulation pattern $P$, we say that $R$ {\it contains} $P$ if some injection from the segments of a representative of $P$ into those of $R$ preserves orientation, incidences, and the relevant order condition (iv).
Otherwise, we say that $R$ {\it avoids} $P$.
Our terminology follows Asinowski–Polley; see \cite[Section 2]{AP25} for background.
For a strong rectangulation pattern $P$, let $\Rec(n; P)$ be the set of strong $n$-rectangulations avoiding pattern $P$.
This paper aims to prove the following theorem.
\begin{thm}\label{thm: upper bound of exponential growth}
    For any strong rectangulation pattern $P$ of size $L\geq 1$, we have
    \begin{equation}\label{eq: upper bound of exponential growth}
        \limsup_{n\to \infty}\# \Rec(n;P)^{1/n} \leq \Lambda - \frac{1}{\Lambda^{3L-1}}.
    \end{equation}
\end{thm}
The following corollary immediately follows.
\begin{cor}\label{cor: proportion is zero}
    For any strong rectangulation pattern $P$ of positive size, we have
    \[
        \lim_{n\to\infty}\frac{\#\Rec(n;P)}{\#\Rec(n)} = 0.
    \]
\end{cor}

Such a pattern avoidance for rectangulations has recently been studied from several complementary angles.
\cite{AP25} initiated a systematic treatment of geometric patterns in rectangulations and completely enumerated classes avoiding the four "$\top$-like" one‑segment patterns and their rotations.
In particular, $\top$-avoiding weak rectangulations are counted by Catalan numbers, while several strong classes are put in bijection with inversion‑sequence classes $I(010,101,120,201)$ and $I(011,201)$.

A second thread characterizes important subclasses by forbidden local configurations ("wall/brick" or "windmill" patterns).
For example, diagonal rectangulations are precisely those that avoid two specific wall patterns, providing a clean avoidance description of this classical class; see Cardinal–Sacrist\'an–Silveira \cite{CSS17}.
In a related direction, \cite{MM23} develops generation frameworks where several guillotine‑type classes are specified by forbidding small sets of local patterns (notably, windmills), by connecting back to permutation encodings.
See also the survey‑style unification in \cite{ACFF24}, where guillotine classes are characterized via avoidance of appropriate mesh patterns under the standard bijections.

These works provide exact enumerations and structural characterizations for a handful of concrete avoidance families.
In contrast, \cref{thm: upper bound of exponential growth} establishes the uniform exponential gap. To our knowledge, no prior result guaranteed such a pattern-uniform exponential deficit for strong rectangulations.
The upper bound in \cref{thm: upper bound of exponential growth} must be non-optimal.
At least, we can improve the exponent $3L-1$ in \cref{eq: upper bound of exponential growth} in many cases.
See \cref{sec: proof of main theorem} for the details.
\begin{que}
    For each strong rectangulation pattern $P$ of size $L$, let
    \[
        \lambda(P)\coloneqq \limsup_{n\to \infty}\#\Rec(n;P)^{1/n}.
    \]
    Can we explicitly determine the values of $\max_{|P|=L}\lambda(P)$ and $\min_{|P|=L}\lambda(P)$?
\end{que}

\subsection*{Organization of this paper}
    In \cref{subsec: notation}, we prepare some notation.
    In \cref{subsec: bijections}, we recall the bijections between the set of rectangulations and the set of leftmost history quadrant excursions for readers' convenience.
    In \cref{sec: proof of main theorem}, we prove \cref{thm: upper bound of exponential growth}.

\subsection*{Acknoledgement}
The author would like to thank Masahiro Nakano for motivating this research, Kengo Nakamura for providing helpful comments to improve the paper, and Shuji Horinaga, Hiroyasu Miyazaki, and Ryoma Onaka for their discussions of the draft.

% ----------------------------------------------------------------

\section{Preliminaries}\label{sec: preliminaries}

% ----------------------------------------------------------------

% ----------------------------------------------------------------

\subsection{Notation}\label{subsec: notation}

% ----------------------------------------------------------------
We adopt the {\it history quadrant walk} conventions of \cite{FNS24}.
See \cite[Section 3]{ACFF24} for a survey treatment.
Let $C = \{ b,r,g,w\}$ be a set of symbols of colors.
A finite sequence
\[
    (v_m)_{m=1}^n = ((h_m,x_m,c_m))_{m=1}^{n} \in \coprod_{n=0}^{\infty}\prod_{m=1}^n(\Z_{\geq 0}^2 \times C)
\]
is a {\it history quadrant walk} (HQW) if for all $m$, the following conditions hold:
\begin{itemize}
    \item $h_m\geq x_m$ for all $m$,
    \item $c_m= b \Rightarrow h_{m+1} = h_m+1$,
    \item $c_m = r$ or $g$ $\Rightarrow h_{m+1} = h_m$, and
    \item $c_m= w \Rightarrow h_{m+1} = h_m - 1$.
\end{itemize}
We regard the null-sequence as a history quadrant walk with the length $n = 0$.
For a history quadrant walk $W = (v_m)_{m=1}^n = ((h_m,x_m,c_m))_{m=1}^{n}$, we say that
$W$ is {\it closed} if $v_n = (0,0,w)$, an {\it excursion} if $W$ is closed and $h_1 = x_1 = 0$, and {\it leftmost} if the following conditions hold:
\[
    \begin{array}{rl}
    c_m \in \{b,r\} \text{ and }c_{m+1}\in \{b,g\} & \Rightarrow x_{m+1} \geq x_m,\\
    \text{otherwise} & \Rightarrow x_{m+1} \geq x_m - 1.
    \end{array}
\]
A history quadrant walk $W = (v_m)_{m=1}^{n}$ is {\it admissible} if $n= 0$, or $v_1\in\{(0,0,b), (0,0,r), (0,0,g)\}$ and $v_n \in \{(0,0,r), (0,0,g), (1,0,w), (1,1,w)\}$. This condition is equivalent to that $W$ is given by excluding the last vertex $(0,0,w)$ from a history quadrant walk.

\begin{rem}
    We use the notation $(h,x,c)$ instead of the original notation $(x,y,c)$ for history quadrant walks. Here, we can convert them to each other by the equality $x + y = h$.
    History quadrant walks can be regarded as a Markov chain with countable vertices $\Z_{\geq 0}^2 \times C$.
    Leftmost history quadrant excursions are paths from $(0,0,c_1)$ to $(0,0,w)$.
    In \cite{ACFF24}, the length of a leftmost history quadrant excursion is defined as the length of the corresponding path.
    In this paper, we call $n$ the length of $W = (v_m)_{m=1}^{n}$.
    Therefore, this length differs by only one from the length handled by \cite{ACFF24}.
\end{rem}

\begin{dfn}
    We use the following notations
    \begin{align}
        \LHQW(n) &\coloneqq \{ \text{leftmost quadrant walk of length } n\},\\
        \LHQWadm(n) &\coloneqq \{ \text{admissible leftmost quadrant walk of length }n\}, \text{ and}\\
        \LHQE(n) &\coloneqq \{ \text{leftmost quadrant excursion of length }n\}.
    \end{align}
    For $E = (v_m)_{m=1}^n \in \LHQE(n)$, let $E^\ast \coloneqq (v_m)_{m=1}^{n-1} \in \LHQWadm(n-1)$.
    For $W=(v_m)_{m=1}^n \in \LHQWadm(n)$, let $\overline{W}$ be the leftmost history quadrant excursion defined as
    \[
        \overline{W} \coloneqq (v_1, \ldots, v_n, (0,0,w)) \in \LHQE(n+1)
    \]
    Then, the operators ${}^\ast$ and ${ }^-$ are bijections between $\LHQW(n+1)$ and $\LHQWadm(n)$.
    For $v \in \Z_{\geq 0}^2 \times C$, $h_0 \in \Z$, and $x_0\in\Z$, let
    \[
        v^{(h_0,x_0)} \coloneqq (h+h_0, x + x_0, c).
    \]
    For $W = (v_m)_{m=1}^{n} \in \LHQW(n)$, let
    \[
        W^{(h_0,x_0)} \coloneqq (v_m^{(h_0,x_0)})_{m=1}^{n}.
    \]
    For $W_0 \in \LHQW(n_0)$ and $X\in \{W,W^\ast, E\}$, we say that
    $W= (v_m)_{m=1}^n \in \LHQX(n)$ avoids $W_0$ if
    \[
        (v_m)_{m=m_0}^{m_0+n_0-1} \neq W_0^{(h,x)}
    \]
    for all $1\leq m_0 \leq n - n_0 + 1$, $h\in \Z$, and $x\in \Z$.
    We set
    \[
        \LHQX(n;W_0) \coloneqq \{W \in \LHQX(n)\ |\ W \text{ avoids }W_0\}.
    \]
    In \cite{ITF09}, a bijection between $\Rec(n)$ and $\LHQE(n)$ is given.
    We denote these bijections by
    \begin{align}
        \proc &\colon \Rec(n) \longrightarrow \LHQE(n)\\
        \pav &\colon \LHQE(n) \longrightarrow \Rec(n)
    \end{align}
    and call them the procedure map and the pavement map.
\end{dfn}

% ----------------------------------------------------------------

\subsection{Procedure map and pavement map}\label{subsec: bijections}

% ----------------------------------------------------------------
In this subsection, we recall the definition of the bijection $\pav$, and explain how it avoids the arbitrariness of the order of the vertices.
See \cite[Section 4.6]{ACFF24} for the details.

\begin{figure}[h]
    \includegraphics[width=4cm]{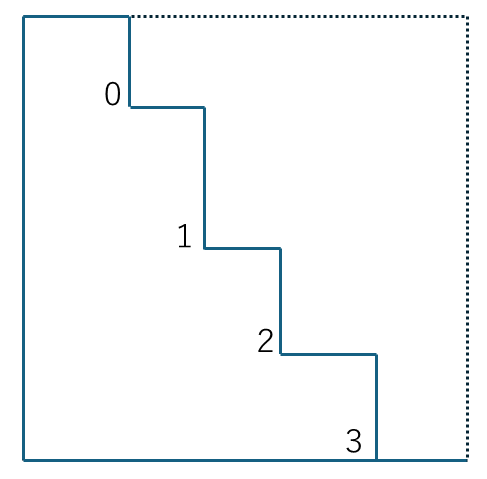}
    \caption{Counting of the steps with $h=3$}
    \label{fig: counting of the steps}
\end{figure}

\begin{figure}[h]
    \includegraphics[width=16cm]{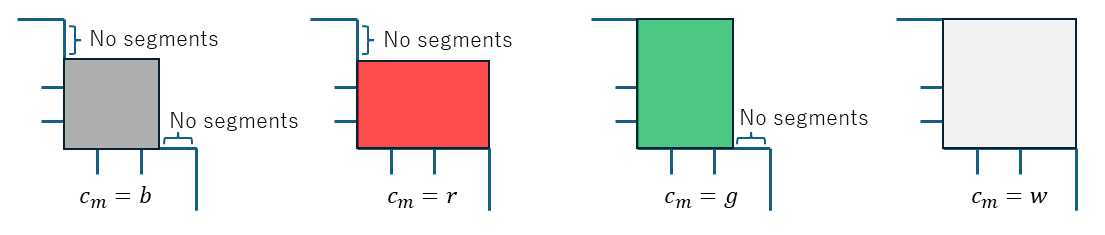}
    \caption{How to place a rectangle for each color}
    \label{fig: how to place}
\end{figure}

For a leftmost history quadrant excursion $E = (v_m)_{m=1}^{n}$, we construct a rectangulation by paving the rectangles corresponding to the excursion steps in the order of index one by one.
In all steps of the pavement, we keep the top and right sides as the shape of staircases.
We count the steps inside the square from $0$ to $h$ as \cref{fig: counting of the steps}.
When the rectangles $v_1, v_2, \ldots, v_{m-1}$ are placed, and $v_m = (h_m,x_m,c_m)$, place the $m$th rectangle at the $x_m$th step of staircases.
As shown in the \cref{fig: how to place}, how to place them is determined by the color $c_m$.
Consequently, we obtain the unique rectangulation corresponding to a given leftmost history quadrant excursion.

\begin{figure}[h]
    \includegraphics[width=4cm]{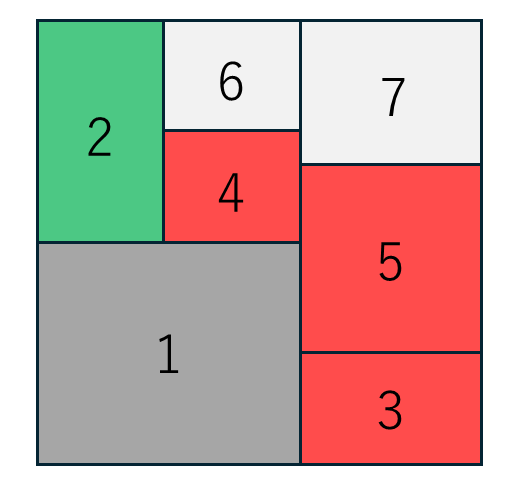}
    \caption{The rectangulation corresponding to $E$}
    \label{fig: rectangulation pattern}
\end{figure}

\cref{fig: rectangulation pattern} shows the rectangulation corresponding to
\begin{equation}\label{eq: E}
    E = \left((0,0,b), (1,0,g), (1,1,r), (1,0,r),(1,1,r), (1,0,w), (0,0,w)\right).
\end{equation}

Conversely, some history quadrant excursions are obtained from one rectangulation as follows.
Consider placing rectangles from the bottom left.
For each rectangle, record the total number of steps of the staircase when placing the rectangle, the number of steps placed, and the color corresponding to the way of placing the rectangle.
To perform this procedure, we have to determine the order of rectangles.
For example, the orders of the rectangles numbered with $2$ and $3$, and $5$ and $6$ in \cref{fig: rectangulation pattern} are not determined by only the word "from the bottom left".
The leftmost property of $E$ forces that the rectangle of $3$ must not be chosen before the rectangle of $2$.
If the rectangle $6$ is placed before the rectangle $5$, we cannot place the rectangle $5$ in any way corresponding to some color.
These observation shows that the order of rectangles such that the corresponding history quadrant walk satisfies the leftmost property is uniquely determined.

% ----------------------------------------------------------------

\section{Proof of the main theorem}\label{sec: proof of main theorem}

% ----------------------------------------------------------------
In this section, we prove \cref{thm: upper bound of exponential growth}.

\begin{proof}
For a strong rectangulation pattern $P$ of size $L \geq 1$, extending segments in one way, we obtain an $(L+1)$-rectangulation $R_0$. Then, we have the following diagram:

\begin{equation}
    \begin{tikzcd}
	{\Rec(n)} & {\LHQE(n)} & {\LHQWadm(n-1)} \\
	{\Rec(n;\seg(R_0))} & {\LHQE(n;\proc(R_0)^\ast)} & {\LHQWadm(n-1; \proc(R_0)^\ast),} \\
	{\Rec(n;P)}
	\arrow["\proc","\sim"', from=1-1, to=1-2]
	\arrow["\ast","\sim"', from=1-2, to=1-3]
	\arrow[hook, from=2-1, to=1-1]
	\arrow[hook, from=2-1, to=2-2]
	\arrow[hook, from=2-2, to=1-2]
	\arrow["\ast","\sim"', from=2-2, to=2-3]
	\arrow[hook, from=2-3, to=1-3]
	\arrow[hook, from=3-1, to=2-1]
\end{tikzcd}
\end{equation}
where we let $\seg(R_0)$ be the strong rectangulation pattern consists of all segments in $R_0$.
Note that since the meaning of pattern avoidance of rectangulations and that of leftmost history quadrant excursions in our sense are different, the map $\Rec(n;\seg(R_0)) \longrightarrow \LHQE(n;\proc(R_0)^{\ast})$ in the diagram is only injective.
To prove the assertion, it is enough to show the inequality
\begin{equation}
    \limsup_{n\to\infty} \left(\# \LHQWadm(n;\proc(R_0)^{\ast})\right)^{1/n} \leq \Lambda - \frac{1}{\Lambda^{3L-1}}.
\end{equation}
If necessary, joining at most $L$ $(0,0,r)$'s to the beginning of $\proc(R_0)^{\ast}$ and at most $L$ $(0,0,g)$'s to the end of $\proc(R_0)^{\ast}$, we obtain an element $W_0 = (v_{0,m})_{m=1}^{L_0}\in \LHQWadm(L_0)$ such that 
\begin{equation}\label{eq: non-coincidence}
    (v_{0,1}, \ldots, v_{0,k}) \neq (v_{0, L_0 - k + 1}^{(h,x)}, \ldots, v_{0,L_0}^{(h,x)})
\end{equation}
for all $h,x \in \Z$ and $1\leq k \leq L_0 - 1$.
\begin{figure}[h]
    \includegraphics[width=10cm]{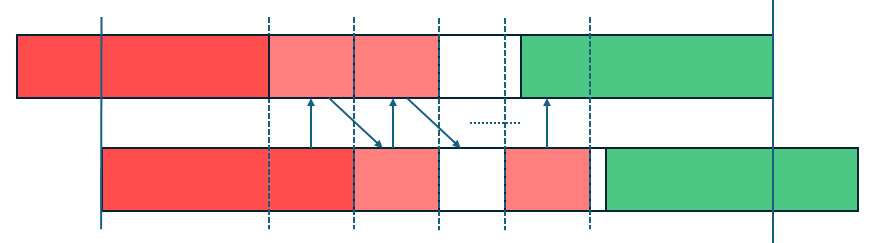}
    \caption{Inheritance of Color in Translation of Walks}
    \label{fig: translation}
\end{figure}
Indeed, if otherwise, all vertices become red as shown in \cref{fig: translation}. However, there are green vertices in our setting.
This is a contradiction.
Note that in many cases, this extension is not needed. Thus, the length $L_0$ can be taken as $L$ in such cases.

Now, it is enough to show the inequality
\begin{equation}
    \limsup_{n\to\infty} \#\LHQWadm(n;W_0)^{1/n} \leq \Lambda - \frac{1}{\Lambda^{L_0 - 1}}.
\end{equation}

For $W = (v_m)_{m=1}^{n} \in \LHQWadm(n; W_0)$ and $q\in \Z_{\geq 0}$, by inserting $q$ translations of $W_0$, we obtain $\binom{n+q}{q}$ distinct elements of $\LHQWadm(n + qL_0)$.
\cref{fig: insertion of a pattern} shows how the insertion of the pattern corresponding to $E^\ast$ defined by \cref{eq: E} looks in the rectangulation.
Let $S(W,q)$ be the set of these elements.
Then, for distinct $W, W' \in \LHQWadm(n; W_0)$, the sets $S(W,q)$ and $S(W',q)$ are disjoint.
Indeed, for an element $\widetilde{W} \in S(W,q)$, we can uniquely reconstruct $W$ by excluding consecutive terms equal to a translation of $W_0$ since we are assuming \cref{eq: non-coincidence}.
\begin{figure}[h]
    \includegraphics[width=0.4\columnwidth]{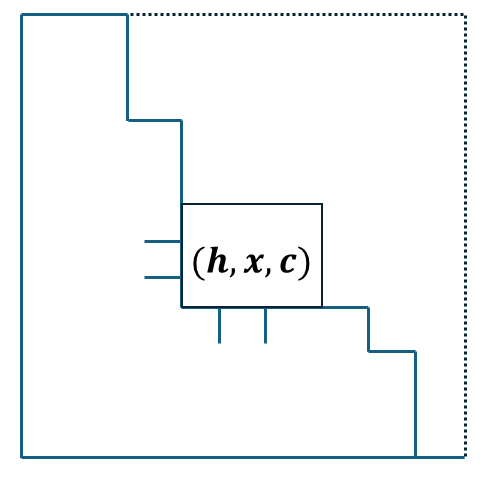}
    \includegraphics[width=0.4\columnwidth]{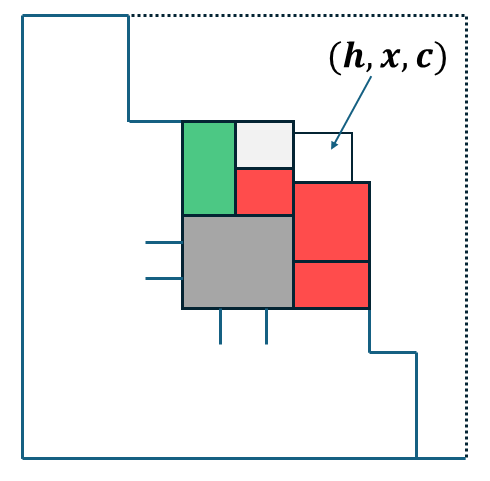}
    \caption{Insertion of a pattern $E^\ast$}
    \label{fig: insertion of a pattern}
\end{figure}
Thus, we have the inequality
\begin{equation}
    \#\LHQWadm(n) \geq \sum_{q=0}^{\floor{n/L_0}}\binom{n-qL_0+q}{q}\cdot \#\LHQWadm(n-qL_0;W_0).
\end{equation}
For each $n\in \Z$, put
\begin{align}
    a_n &= \# \LHQWadm(n),\ \text{and}\\
    b_n &= \# \LHQWadm(n;W_0),
\end{align}
where we regard $a_n = b_n = 0$ for negative $n$, and $a_0 = b_0 = 1$.
Consider the formal series
\begin{align}
    A(z) &= \sum_{n=0}^{\infty} a_n z^n,\\
    B(z) &= \sum_{n=0}^{\infty} b_n z^n.
\end{align}
Then, we already know that the convergence radius of $A(z)$ is $1/\Lambda$ by \cref{thm: asymptotics of rec n}.
By the equality
\[
    \sum_{q\geq 0} \binom{m+q}{q} t^q = (1-t)^{-(m+1)},
\]
we obtain the inequalities
\begin{align}
    A(z) &\succeq \sum_{n=0}^{\infty} \left[ \sum_{q=0}^{\floor{n/L_0}} \binom{n - qL_0 + q}{q} b_{n-qL_0}\right] z^n\\
    &= \sum_{m=0}^{\infty} \sum_{q=0}^{\infty} \binom{m+q}{q} z^{qL_0}b_mz^m &\text{by putting }m = n-qL_0\\
    &= \sum_{m=0}^{\infty} (1-z^{L_0})^{-(m+1)} b_m z^m\\
    &= \frac{1}{1-z^{L_0}}\cdot B\left(\frac{z}{1-z^{L_0}}\right).
\end{align}
Since $A(z)$ is holomorphic on $|z|<1/\Lambda$ by \cref{thm: asymptotics of rec n}, so is the right-hand side.
Thus, $B(z)$ is holomorphic on $|z|<r$ for 
\[
    r = \sup\left\{ \frac{z}{1-z^{L_0}}\ | \ |z|<\frac{1}{\Lambda}\right\}
    = \frac{\frac{1}{\Lambda}}{1 - \frac{1}{\Lambda^{L_0}} }.
\]
Consequently, we obtain the inequality
\[
    \limsup_{n\to\infty} b_n^{1/n} \leq 1/r = \Lambda - \frac{1}{\Lambda^{L_0 - 1}}.
\]

\end{proof}

% --------------------------------------------------------------------------
%		References
% --------------------------------------------------------------------------
%\newpage
\bibliographystyle{alpha}
\bibliography{rectangulation}

% --------------------------------------------------------------------------
\end{document}